\documentclass[a4paper,12pt]{article}

\usepackage{hyperref}

\usepackage[top=3cm,left=2.5cm,right=2.5cm,bottom=3cm]{geometry}
\usepackage{relsize}
\usepackage{lineno}

\usepackage{amsmath,amsthm,pb-diagram,amssymb,comment}
\usepackage{amsfonts,graphicx,color}
\usepackage{enumitem, fancyhdr, dsfont, multicol}
\usepackage{thmtools,cleveref}
\usepackage[normalem]{ulem}
\usepackage{stmaryrd}
\usepackage{textcomp}
\usepackage{mathtools}
\usepackage{fontawesome}
\usepackage{textcomp}

\usepackage{tikz,float}
\usepackage{mleftright}
\usepackage{mathrsfs}
\usepackage{tikz-cd}
\usepackage{dsfont}
\usepackage{stmaryrd}

\hypersetup{
    colorlinks=true,
    linkcolor=magenta,
    filecolor=magenta,
    urlcolor=magenta,
    citecolor=magenta,
}
\definecolor{caribbeangreen}{rgb}{0.0, 0.8, 0.6}

\setlist{topsep=0ex,itemsep=1ex}

    \DeclareMathOperator{\dom}{{\rm dom}}

    \DeclareMathOperator{\ran}{{\rm ran}}

    \newcommand{\Swf}{\mathcal{S}}

    \newcommand{\Ba}{\mathcal{B}\boldsymbol{a}}

    \newcommand{\RRR}{\mathds R}
\newcommand{\QQQ}{\mathds Q}
\newcommand{\ZZZ}{\mathds Z}

    \newcommand{\menos}{\smallsetminus}
    
    \DeclareMathOperator{\pts}{\mathcal{P}}

    \DeclareMathOperator{\cov}{{\rm cov}}

    \DeclareMathOperator{\cf}{{\rm cf}}

    \newcommand{\la}{\langle}
    \newcommand{\ra}{\rangle}

\newcommand{\Fr}{\mathrm{Fr}}

\newcommand{\Lc}{\mathbf{Lc}}

\newcommand{\id}{\mathrm{id}}

\newcommand{\baire}{{}^{\omega} \omega}

\newcommand{\set}[2]{\left\{#1 \colon #2\right\}}

\newcommand{\Fn}{\mathrm{Fn}}

\newcommand{\cantor}{{}^\omega2}

\newcommand{\FAM}{\mathrm{FAM}}

\newcommand{\alt}{\mathrm{ht}}

\newcommand{\Int}{\mathrm{int}}

\newcommand{\Ord}{\mathrm{Ord}}

\newcommand{\Leb}{\mathrm{Leb}}

\newcommand{\forces}{\Vdash}

\newcommand{\Lev}{\mathrm{Lev}}

\newcommand{\suc}{\mathrm{succ}}

\newcommand{\bfm}{\mathbf{m}}

\newcommand{\bbD}{\mathbb{D}}
\newcommand{\bbE}{\mathbb{E}}

\newcommand{\bbP}{\mathbb{P}}
\newcommand{\bbQ}{\mathbb{Q}}
\newcommand{\bbR}{\mathbb{R}}

\newcommand{\cB}{\mathscr{B}}
\newcommand{\cC}{\mathscr{C}}

\newcommand{\gb}{\mathfrak{b}}

\newcommand{\calT}{\mathcal{T}}

\newcommand{\calA}{\mathcal{A}}

\newcommand{\calB}{\mathcal{B}}
\newcommand{\calP}{\mathcal{P}}
\newcommand{\calS}{\mathcal{S}}
\newcommand{\calN}{\mathcal{N}}

\newcommand{\calX}{\mathcal{X}}

\newcommand{\finseqb}{{}^{ < \omega } \omega}
\newcommand{\finseq}{{}^{ < \omega } 2}
\newcommand{\finseqz}{{}^{ < \omega } Z}

\newcommand{\varp}{\varepsilon}

\newcommand{\rest}{\restriction}

\newenvironment{PROOF}[2][\proofname.]
   {\begin{proof}[#1]\renewcommand{\qedsymbol}{\textsquare$_{\mathrm{#2}}$}}
   {\end{proof}}

\renewcommand{\setminus}{\smallsetminus}
\newcommand{\conj}{\wedge}

\DeclareMathOperator{\loss}{loss}
\DeclareMathOperator{\trunk}{trunk}

\definecolor{sub0}{RGB}{29,32,137}
\definecolor{sub1}{RGB}{1,71,157}
\definecolor{sub2}{RGB}{1,104,183}
\definecolor{sub3}{RGB}{0,160,234}
\definecolor{sug}{RGB}{0,154,68}
\definecolor{suy}{RGB}{208,219,1}

\definecolor{redun}{rgb}{0.65, 0.11, 0.19}
\definecolor{greenun}{rgb}{0.58, 0.71, 0.23}
\definecolor{dodger}{rgb}{0.0,0.5,1.0}
\definecolor{carrotorange}{rgb}{0.93, 0.57, 0.13}

\begin{document}

\makeatletter
\def\@roman#1{\romannumeral #1}
\makeatother

\newcounter{enuAlph}
\renewcommand{\theenuAlph}{\Alph{enuAlph}}

\theoremstyle{plain}
  \newtheorem{theorem}{Theorem}[section]
  \newtheorem{corollary}[theorem]{Corollary}
  \newtheorem{lemma}[theorem]{Lemma}
  \newtheorem{mainlemma}[theorem]{Main Lemma}
  \newtheorem{mainproblem}[theorem]{Main Problem}
  \newtheorem{construction}[theorem]{Construction}
  \newtheorem{prop}[theorem]{Proposition}
  \newtheorem{clm}[theorem]{Claim}
  \newtheorem{fact}[theorem]{Fact}
  \newtheorem{exer}[theorem]{Exercise}
  \newtheorem{question}[theorem]{Question}
  \newtheorem{problem}[theorem]{Problem}
  \newtheorem{cruciallem}[theorem]{Crucial Lemma}
  \newtheorem{conjecture}[theorem]{Conjecture}
  \newtheorem{assumption}[theorem]{Assumption}
  \newtheorem{teorema}[enuAlph]{Theorem}
  
  \newtheorem*{thm}{Theorem}
    \newtheorem*{thma*}{Theorem A}
    \newtheorem*{thmb*}{Theorem A}
    \newtheorem*{thmc*}{Theorem B}
    \newtheorem*{thmd*}{Theorem D}
 
  \newtheorem*{corolario}{Corollary}
\theoremstyle{definition}
  \newtheorem{definition}[theorem]{Definition}
  \newtheorem{example}[theorem]{Example}
  \newtheorem{remark}[theorem]{Remark}
    \newtheorem{hremark}[theorem]{Historical Remark}
    \newtheorem{observation}[theorem]{Observation}
    \newtheorem{notation}[theorem]{Notation}
  \newtheorem{context}[theorem]{Context}

\parindent 0pt

  \newtheorem*{defi}{Definition}
  \newtheorem*{acknowledgements}{Acknowledgements}

\numberwithin{equation}{theorem}
\renewcommand{\theequation}{\thetheorem.\arabic{equation}}

\def\sectionautorefname{Section}
\def\subsectionautorefname{Subsection}


\title{The intersection number for forcing notions}


\author{\Large Andrés F. Uribe-Zapata\\
\small \href{mailto:andres.zapata@tuwien.ac.at}{andres.zapata@tuwien.ac.at}\\ 
\small \href{https://sites.google.com/view/andres-uribe-afuz/home}{https://sites.google.com/view/andres-uribe-afuz/home}
}

\date{Institute of Discrete Mathematics and Geometry\\
Faculty of Mathematics and Geoinformation\\
TU Wien\\ 
 Wiedner Hauptstrasse 8--10, A--1040\\
Vienna, Austria\\ \vspace{0.6cm}
\small \today}


\maketitle

\begin{abstract}
    \noindent
    Based on works of Saharon Shelah, Jakob Kellner, and Anda T\u{a}nasie for controlling the cardinal characteristics of the continuum in ccc forcing extensions, in the author's master's thesis was introduced a new combinatorial notion: the \emph{intersection number for forcing notions}, which was used in such thesis to build a general theory of iterated forcing using finitely additive measures. In this paper, we present the definition of such a notion and prove some of its fundamental properties in detail. Additionally, we introduce a new linkedness property called \emph{$\mu$-intersection-linked}, prove some of its basic properties, and provide some interesting examples.

    \vspace{0.37cm}

\end{abstract}

\section{Introduction}

The notion of \emph{intersection number} appears in different contexts in mathematics, for example, in algebraic geometry\footnote{Where there is in fact a sub-branch called \emph{intersection theory} (see \cite{Fulton}).}, graph theory (see e.g. \cite{Graph}), and mathematical logic. This notion has even been taken into theoretical physics contexts, where some interesting applications have appeared (see \cite{Physic}). In mathematical logic, particularly in the context of Boolean algebras, such a concept first appeared in John Kelley's work \cite{Kelley59} to deal with the general problem of the existence of finitely additive non-trivial measures on Boolean algebras. There, motivated by a conjecture of Alfred Tarski, John Kelley succeeded in characterizing the existence of finitely additive measures using the notion of intersection number: he proved that, given a Boolean algebra $\cB$, there exists a finite strictly positive finitely additive measure on $\cB$ if, and only if, $\cB^{+} \coloneqq \cB \setminus \{ 0_{\cB} \}$ can be decomposed as a union of countable many subsets with positive intersection number (see~\cite[Thm.~4]{Kelley59}). Much later, in 2023, motivated by the ideas of John Kelley, in the author's master's thesis the intersection number appeared in a new context: forcing theory. In this case, the so-called \emph{intersection number for forcing notions} played a fundamental role in presenting a general theory of iterated forcing using finitely additive measures in the following way. In 2000, Saharon Shelah built a forcing iteration using finitely additive measures to show that, consistently, the covering number of the Lebesgue null ideal, $\cov(\calN),$ may have countable cofinality. The iteration used random forcing, whose structure was essential both to be able to extend the iteration in each step and to apply it to the problem of the cofinality of $\cov(\calN).$ Later, in 2019, Jakob Kellner, Saharon Shelah, and Anda T\u anasie \cite{KST} made new contributions to the method of iterations using finitely additive measures. In particular, they defined a new notion, called \emph{strong FAM limit for intervals} (see \cite[Def.~1.7]{KST}), which turned out to be the key to extending iterations using finitely additive measures at successor steps (see \cite[Lem.~2.25]{KST} and \cite[Thm.~4.3.16]{uribethesis}). Finally, in 2023, in the author's master's thesis (see \cite[Def.~4.1.1]{uribethesis}), it was shown that the intersection number was the key to extending iterations using finitely additive measures at limit steps (see \cite[Thm.~4.3.18]{uribethesis}). Based on the ideas of the intersection number for forcing notions and the strong-$\FAM$-limits, a new linkedness property called \emph{$\theta$-$\FAM$-linkedness} was defined in \cite[Def.~4.2.8]{uribethesis}, where $\theta$ is an infinite cardinal, to present a general theory of iterated forcing using finitely additive measures (see \cite[Ch.~4]{uribethesis}). In particular, due to results about the intersection number, this property allows to extend the iterations in the limit and successor steps, preserving the unbounded number $\gb$, and preserving what in \cite{CMU} is defined as \emph{strongly-$\theta$-anti-Bendixson-families}, which is essential to prove the consistency of $\cf(\cov(\calN)) = \aleph_{0}$.

Currently, together with Miguel A. Cardona, there is an article in preparation \cite{CMU} based on \cite{uribethesis}, where a general theory of iterated forcing using finitely additive measures is presented. As mentioned before, the intersection number for forcing notions appears there, and to build the aforementioned theory it is important to have available properties and results about it, especially for the proof of the extension theorem at limit steps (see~\cite[Thm.~4.3.18]{uribethesis}). However, due to length constraints, the authors in \cite{CMU} decided to omit the proofs of most of these results, and this is where the central idea of this paper was born, which can be thought of as a detailed  presentation of the proofs omitted there. This paper is structured as follows. In \autoref{2}, we present the necessary preliminary notions for the development of the subsequent sections, in particular, we present the definition of the forcing notion $\tilde{\bbE}$ from \cite{KST}. In \autoref{3}, all the details about the intersection number that are omitted in \cite{CMU} are presented. In particular, it is shown how the intersection number for forcing notions is related to the intersection number for Boolean algebras through forcing-completions (see \autoref{s3}), and the behavior of the intersection number under complete embeddings is studied (see \autoref{s6}). In addition, some basic properties of the intersection number are proved. In \autoref{4}, a new linkedness property is introduced: \emph{$\theta$-intersection-linked}, which arises naturally by weakening the definition of $\theta$-$\FAM$-linked without considering the property referring to the $\FAM$-limits. In particular, based on \cite{MU23}, we show that any Boolean algebra with a strictly positive probability measure with the $\theta$-density property is $\theta$-intersection-linked. As a consequence, the measure algebra adding $\theta$-many random reals is $\theta$-intersection-linked.

\section{Preliminaries}\label{2}

We denote by $\ZZZ,$ $\QQQ$, and $\RRR$ the sets of integers, rational and real numbers respectively. Notice the difference between the symbols ``$\QQQ$'' and ``$\bbQ$'', and ``$\RRR$'' and ``$\bbR$'', since $\bbQ$ and $\bbR$ will be used to denote forcing notions. If $I \subseteq \RRR$ is an interval, we define $I_{\QQQ} \coloneqq I \cap \QQQ.$ We use the symbol $\omega$ to denote the set of natural numbers, $\Ord$ denotes the proper class of ordinal numbers, and $\mathrm{Card}$ denotes the proper class of the cardinal numbers.

Let $A, B$ be sets, $\alpha \in \Ord$ and $\kappa, \lambda \in \mathrm{Card}.$ Define $$[A]^{\kappa} \coloneqq \{ X \subseteq A \colon \vert X \vert = \kappa\} \text{ and } [A]^{< \kappa} \coloneqq \{ X \subseteq A \colon \vert X \vert < \kappa \}.$$  We denote by ${}^{A} B$ the set of functions $f$ from $A$ into $B$ and ${}^{< \alpha} A \coloneqq \bigcup \{ {}^{\xi} A \colon \xi < \alpha  \}.$ Similarly, ${}^{\leq \alpha}A \coloneqq {}^{< \alpha +1 }A.$ For any $t \in {}^{< \alpha} A$ we define its \emph{length} by $\lg(t) \coloneqq \dom(t).$  We use the symbols ``$\langle$'' and ``$\rangle$'' to denote sequences and ``$\langle \ \rangle$'' to denote the empty sequence. If $A, B$ are non-empty, then $\Fn(A, B)$ is the set of \emph{finite partial functions} from $A$ into $B,$ that is, functions $f \colon X \to B,$ such that $X \in [A]^{< \aleph_{0}}.$ If $\la A, R\ra$  is a preoder and $B \subseteq A,$ then we define $B^{\uparrow} \coloneqq \{ a \in A \colon \exists b \in B\ (b \mathrel{R} a) \}.$ Again, if $B \subseteq A,$  we denote by $\chi_{B}$ the characteristic function of $B$ in $A.$ 

\subsection{Finitely additive measures on Boolean algebras}\label{2.5} 

In this subsection, based on \cite[Ch.~3]{uribethesis} and  \cite{CMUP}, we present some basic notions about finitely additive measures on Boolean algebras. 

Let $\cB$ be a Boolean algebra. Recall that, if $I$ is an ideal on $\cB,$ then we can define the quotient between $\cB$ and $I,$ which we denote by $\cB / I,$ and the equivalent classes in $\cB / I$ are denoted by $[b]_{I}$ for any $b \in \cB.$  

\begin{definition}
    Let $\cB$ be a Boolean algebra. A \emph{finitely additive measure on $\cB$} is a function $\Xi \colon \cB \to [0,\infty]$ satisfying the following conditions:   
    
\begin{enumerate}[label=(\arabic*)]
    \item \label{m4a} $\Xi (0_{\cB})=0$,
            
    \item \label{m4b} $\Xi(a\vee b)=\Xi(a)+\Xi(b)$ whenever $a,b\in\cB$ and $a \wedge b= 0_{\cB}$.
\end{enumerate}

    We say that $\Xi$ is a \emph{measure on $\cB$} if it satisfies~\ref{m4a} and
    \begin{enumerate}
        \item[(2)'] $\Xi \! \left(\bigvee\limits_{n<\omega} b_n\right)=\sum\limits_{n<\omega} \Xi(b_n)$ whenever $\{ b_n \colon n < \omega \}  \subseteq \cB$, $\bigvee_{n<\omega}b_n$ exists, and $b_m\wedge b_n=0_\cB$ for $m\neq n$.
    \end{enumerate}
        
    We exclude the trivial finitely additive measure, that is, we will always assume $\Xi(1_\cB)>0$. If $\Xi(1_{\cB}) < \infty,$ we say that $\Xi$ is \emph{finite}.  When $\Xi(b)>0$ for any $b \in \cB^{+}$ we say that $\Xi$ is \emph{strictly positive}.  Also, if $\Xi(1_\cB)=1,$ we say that $\Xi$ is a \emph{finitely additive probability measure}.
\end{definition}

\begin{definition}
    Let $\cB$ be a Boolean algebra, $\Xi$ a finitely additive measure on it, and $b \in \cB$ with positive finite $\Xi$-measure. We define the function $\Xi_{b} \colon \cB \to [0, \infty]$ by $\Xi_{b}(a) \coloneqq \frac{\Xi(a \wedge b)}{\Xi(b)}$ for any $a \in \cB$.
\end{definition}

It is clear that $\Xi_{b}$ is a finitely additive probability measure.

Based on \cite[Def.~5.4]{BCM2}, we introduce the notion of \emph{$\theta$-density property}: 


\begin{definition}\label{m14}
    Let $\cB$ a Boolean algebra, $\mu$ a strictly positive probability measure on $\cB,$ and $\theta$ an infinite cardinal. We say that $\mu$ satisfies the \emph{$\theta$-density property} if there exists an $S \subseteq \cB^{+}$ such that $\vert S \vert \leq \theta$ and, for any $\varp > 0$ and $b \in \cB^{+}$, there is some $s \in S$ such that $\mu_{s}(b) \geq 1 - \varp.$
\end{definition}

For example, any $\sigma$-centered Boolean algebra has a strictly positive finitely additive measure with the $\aleph_{0}$-density property (see \cite[Lem.~5.5]{BCM2}).

It is clear that we have the following hereditary property: 

\begin{lemma}\label{m15.1}
    Let $\cB,$ $\cC$ be Boolean algebras, $\mu$ a strictly positive finitely additive probability measure on $\cB$ and $\theta$ an infinite cardinal. Assume that $\cC$ is a Boolean subalgebra of $\cB$ and $\mu$ has the $\theta$-density property witnessed by $S.$ If $S \subseteq \cC$, then $\mu \rest \cC$ also has the $\theta$-density property witnessed by $S.$
\end{lemma}

Recall that a \emph{field of sets over $X$} is a Boolean subalgebra of $\pts(X)$, where the latter is a Boolean algebra with the set-theoretic operations. Now, we introduce integration over a field of sets. Fix a non-empty set $X$, a field of sets $\cB$ over $X$, and a finitely additive measure $\Xi \colon \cB \to [0, \infty).$ Motivated by the definition of Riemann's integral, if  $f \colon X \to \mathbb{R}$ is a bounded function we can naturally define $\int_{X} f d \Xi$ if it exists. In this case, we say that $f$ is \emph{$\Xi$-integrable}  (see \cite[Def.~3.5.3]{uribethesis}).  For example, any bounded function is $\Xi$-integrable when $\dom(\Xi) = \calP(X)$ (see \cite[Thm.~3.5.10]{uribethesis}). In fundamental aspects, the integral with respect to finitely additive measures behaves similarly to the Riemann integral, that is, we have available the basic properties of the integral such as those presented in \cite[Sec.~3.5]{uribethesis}. Next, we list some of these properties which are necessary for the development of this work.

\begin{lemma}\label{t47}
    Let $f, g$ be $\Xi$-integrable functions and $c \in \bbR.$ Then 

    \begin{enumerate}[label=\rm{(\arabic*)}]
        \item  $fg$ and $cf$ are $\Xi$-integrable. In this case,  $\int_{X} (cf) d \Xi = c \int_{X} f d \Xi.$

        \item \label{t50} Let $\{ f_{i} \colon i < n \}$ a finite sequence of $\Xi$-integrable functions. Then $\sum_{i < n} f_{i}$ is $\Xi$-integrable and $$\int_{X} \left( \sum_{i < n} f_{i}\right) d \Xi = \sum_{i < n} \left( \int_{X} f_{i} \, d \Xi \right) \!.$$

        \item \label{t53} If $f, g \colon X \to \RRR$ are $\Xi$-integrable functions and $f \leq g,$ then $ \int_{X} f d \Xi \leq \int_{X} g d \Xi.$

        \item \label{t63} If $E \in \cB,$ then $\chi_{E}$ is $\Xi$-integrable and $\int_{X} \chi_{E}  d \Xi = \Xi(E).$
    \end{enumerate}
\end{lemma}

\subsection{Trees}

Let $Z$ be a non-empty set. A \emph{subtree of $\finseqz$} is a set $\calT \subseteq \finseqz$ such that, for any $\rho, \eta \in \finseqz,$ if $ \rho \subseteq \eta$ and $\eta \in \calT,$ then $\rho \in \calT.$ For instance, for any $n^{\ast} < \omega, \, {}^{n^{\ast} \geq }2$ is a subtree of $\finseq$ called the \emph{complete binary tree} of height $n^{\ast} +1.$  

Fix a tree $\calT$ and $\rho, \eta \in \calT.$ Elements in $\calT$ are called \emph{nodes}, and if $\rho \subseteq \eta,$ we say that $\eta$ is an \emph{extension of $\rho$}. We say that $\rho$ and $\eta$ are \emph{compatible} if they are compatible as functions.  The \emph{height of $\rho$ in $\calT$} is $\alt_{\calT}(\rho) \coloneqq \dom(\rho)$ and  the \emph{height of $\calT$}  is $\alt(\calT) \coloneqq \sup \{ \alt_{\calT}(\rho) + 1 \colon \rho \in \calT \}.$ When the context is clear, we  simply denote $\alt_{\calT}(\rho)$ as $\alt(\rho)$. An \emph{infinite branch} of $\calT$ is an element of $z \in {}^{\omega}Z$ such that, for any $n < \omega, \, z \, {\rest} \, n \in \calT.$ The set of infinite branches of $\calT$ is denoted by $[\calT].$ Also, for any $\rho \in \finseqz,$  we define $[\rho] \coloneqq \{x  \in {}^{\omega} Z \colon \rho \subseteq  x \}. $  

Below we define other notions related to tress: 

\begin{definition}\label{y8}

    Let $\calT$ be a subtree of $\finseqz$ and $\rho \in \calT.$ We define:

    \begin{enumerate}[label=\rm{(\arabic*)}]
        \item For any $h < \omega,$ the \emph{$h$-th level of $\calT$} is  $\Lev_{h}(\calT) \coloneqq \calT \cap {}^{h} Z,$ 

        \item $\calT_{\geq \rho} \coloneqq \{ \eta \in \calT \colon \rho \subseteq \eta \}$ is the set of \emph{successors of $\rho$} in $\calT.$ 

        \item $\suc_{\calT}(\rho) \coloneqq  \calT_{\geq \rho} \cap \Lev_{\alt_{\calT}(\rho) + 1}(\calT),$ that is, $\suc_{\calT}(\rho)$ is the set of \emph{immediate successors of $\rho$} in $\calT.$ 

        \item $\max(\calT) \coloneqq \{ \rho \in \calT \colon \suc_{\calT}(\rho) = \emptyset \},$ that is, it is the set of \emph{maximal nodes} of $\calT.$ 


    \end{enumerate}

    When the context is clear we simply write  ``$\suc(\rho)$'' instead of ``$\suc_{\calT}(\rho)$''.
\end{definition}

Notice that, for any $\rho \in \calT, \, \alt(\rho) < \omega$ and, if $[\calT] \neq \emptyset,$ then $\alt(\calT) = \omega.$  Also, $\rho \in \Lev_{j}(\calT)$ if, and only if, $\alt(\rho) = j.$

\subsubsection{Forcing notions} 

We assume the reader to be familiar with basic techniques of set theory (see e.g.~\cite{Kunen}), in particular with forcing theory, however, we review some notation and definitions.

\begin{definition}\label{n4}
    Let $(\bbP, \leq)$ a forcing notion. The \emph{separative order induced by $\leq$ on $\bbP$}, denoted by $\leq^{\bullet},$ is defined by: $p \leq^{\bullet} q $ if, and only if, for any $r \leq p, \,  r \parallel_{\bbP} q,$ that is, $r$ and $q$ are \emph{compatible} in $\bbP.$ 
\end{definition}

Notice that, if $p \leq q,$ then $p \leq^{\bullet} q,$ however, the converse is not true in general. When the converse holds, we say that $\bbP$ is a \emph{separative} forcing notion. For example, Boolean algebras are separative. A converse attempt at forcing notions, in general, can be the following lemma: 

\begin{lemma}\label{n5}
    Let $\bbP$ be a forcing notion, $p \in \bbP$ and $\langle p_{i} \colon i < n \rangle \subseteq \bbP.$ If for every $i < n, \, p \leq^{\bullet} p_{i},$ then there exists some $q \leq p$ such that, for any $i < n, \, q \leq p_{i}.$  
\end{lemma}


Diego A. Mej\'ia  \cite{mejiavert} introduced  the notion of \emph{Fr\'echet-linked}:

\begin{definition}\label{n73}
    Let $\bbP$ be a forcing notion and $\mu$ an infinite cardinal.
    
    \begin{enumerate}[label=\rm{(\arabic*)}]
        \item We say that a set $Q \subseteq \bbP$ is \emph{Fr\'echet-linked in $\bbP$}, abbreviated \emph{$\Fr$-linked} if, for any sequence $\bar{p} = \langle p_{n} \colon n < \omega \rangle \in {}^{\omega} Q,$ there exists some condition $q \in \bbP$ such that $q \forces_{\bbP} \text{``} \vert \{ n < \omega \colon p_{n} \in \dot{G} \} \vert = \aleph_{0}$''. 

        \item We say that $\bbP$ is \emph{$\mu$-$\rm{Frechet}$-linked}, abbreviated by \emph{$\mu$-$\Fr$-linked}, if $\bigcup_{\alpha < \mu} Q_{\alpha}$ is dense in $\bbP$ for some sequence   $\langle Q_{\alpha} \colon \alpha < \mu \rangle$  of $\Fr$-linked subsets of $\bbP.$
    \end{enumerate}
\end{definition}

Denote the ground model by $V$. A \emph{dominating real over $V$} is a real number $d \in \baire$ such that, for any $x\in\omega^\omega\cap V,$ $x$ is eventually dominated by $d$. The $\Fr$-linkedness notion arises implicitly in Arnold Miller's proof (see \cite{Miller2}) that the \emph{eventually different forcing notion} $\bbE$ does not add dominating reals, which turned out to be very useful because it implies preservation properties (see \cite[Sec.~3]{mejiavert}).

\begin{theorem}\label{n80}
    If $\bbP$ is $\sigma$-$\Fr$-linked, then $\bbP$ does not add dominating reals. 
\end{theorem}

Recall that random forcing is the complete Boolean algebra $\calB(\cantor) / \calN(\cantor),$ where $\calB(\cantor)$ is the Borel $\sigma$-algebra on $\cantor$ and $\calN \coloneqq \calN(\cantor)$ is the ideal of Lebesgue-null sets in $\cantor.$ More generally, Kenneth Kunen in~\cite{Ku84} presents an approach to deal with measure and category in the context of ${}^{I}2$ with $I$ possibly uncountable. With the product topology of the discrete space $\{ 0, 1 \}$, the basis is described by the basic clopen sets $[s] \coloneqq \{x \in {}^{I} 2 \colon  s \subseteq x \}$ for $s \in \Fn(I,2)$, and $\Ba({}^{I}2)$, the $\sigma$-algebra of \emph{Baire sets}, is defined as the $\sigma$-algebra generated by the clopen sets (which coincides with the product $\sigma$-algebra). We can also define the product measure $\Leb^{I}$ on the $\sigma$-algebra $\Ba({}^{I}2)$, where $\{ 0, 1 \}$ is endowed with the measure whose value on $\{ 0 \}$  and $\{ 1\}$  is $\frac{1}{2}$.  Then $\calN({}^{I} 2)$, the set of $\Leb^{I}$-null sets in ${}^{I} 2$, is an ideal on $\Ba({}^{I} 2)$ and, based on \cite[Def.~IV.7.33]{Kunen}, we can define: 

\begin{definition}\label{i13}
    For any infinite cardinal $\theta,$ the \emph{measure algebra adding $\theta$-many random reals} is the quotient $\cB_{\theta} \coloneqq  \Ba({}^{\theta} 2) / \calN({}^{\theta} 2).$ 
\end{definition} 

If, abusing the notation, for any infinite cardinal $\theta,$ we define $\Leb^{\theta} \colon \cB_{\theta} \to [0, \infty]$ such that $\Leb^{\theta}([A]_{\calN({}^{\theta} 2)}) \coloneqq \Leb^{\theta}(A)$, then   

\begin{theorem}\label{i3.0}
    For any infinite cardinal $\theta,$ $\Leb^{\theta}$ is a strictly positive probability measure on $\cB_{\theta}$ satisfying the $\theta$-density property. 
\end{theorem}

A proof of this theorem can be found in \cite{MU23}.  

Below, we introduce a variant of the eventually different forcing notion $\bbE$ introduced in~\cite[Def.~1.12]{KST} as a variation of a forcing notion by Haim Horowitz and Saharon Shelah (see \cite{HSh}). 

\begin{definition}\label{v48}
    By induction on the height $h \geq 0,$ we define a countable tree $\tilde{\calT} \subseteq \finseqb,$ functions $\varrho, \pi, a, M \colon \omega \to \omega$ and a map $\mu_{h} \colon M(h)+1 \to \RRR$ as follows: 
    \begin{enumerate}
        \item $\Lev_{0}(\tilde{\calT}) \coloneqq \{ 
        \langle \ \rangle \},$ that is, the unique element of height $0$ is the empty sequence $\langle \ \rangle,$ and: 

        \begin{multicols}{2} 
            \begin{itemize}
                \item $\varrho(0) \coloneqq 2,$
        
                \item $\pi(0) \coloneqq 2,$
        
                \item $a(0) \coloneqq 4,$

                \item $M(0) \coloneqq 16,$
        
                \item $\mu_{0}(n) \coloneqq \log_{4} \left(  \frac{16}{16-n} \right)$ for $n < 16,$
        
                \item $\mu_{0}(16) \coloneqq \infty.$
            \end{itemize}
        \end{multicols}

        \item Assume we have defined $\Lev_{h}(\tilde{\calT})$ for $h < \omega.$ For any $\rho \in \Lev_{h}(\tilde{\calT}),$ define $$ \suc_{\calT}(\rho) \coloneqq \{ \rho^{\smallfrown} \langle \ell \rangle \colon \ell \in M(h) \} $$ and $$ \Lev_{h+1}(\tilde{\calT}) \coloneqq \bigcup_{\rho \in \Lev_{h}(\tilde{\calT})} \suc_{\tilde{\calT}}(\rho). $$ Now, 

        \begin{itemize}
            \item $\varrho(h+1) \coloneqq \max \{ \vert \Lev_{h+1}(\tilde{\calT} )\vert, \, h +3 \},$

            \item $\pi(h+1) \coloneqq [(h+2)^{2} \varrho(h+1)^{h+2}]^{\varrho(h+1)^{h+1}},$

            \item $a(h+1) \coloneqq \pi(h+1)^{h+3},$

            \item $M(h+1) \coloneqq a(h+1)^{2},$ 

            \item $\mu_{h+1}(n) \coloneqq \left\{ \begin{array}{lcc}
            \log_{a(h+1)}  \left(  \frac{M(h+1)}{M(h+1) - n}\right), &  \rm{if}  & 0 \leq n < M(h+1), \\
            \\ \infty, &  \rm{if}  & n = M(h+1).
            \end{array}
            \right.$

            \item If $\rho \in \Lev_{h+1}(\tilde{\calT})$ and $A \subseteq \suc({\rho}),$ then we set $\mu_{\rho}(A) \coloneqq \mu_{n+1}(\vert A \vert).$ 
        \end{itemize}
    \end{enumerate}

    Finally, define $$\tilde{\calT} \coloneqq \bigcup_{h < \omega} \Lev_{h}(\tilde{\calT}),$$  and if $p \subseteq \tilde{\calT}$ is a tree and $\rho \in p,$ $\mu_{\rho}(p) \coloneqq \mu_{\rho}( \suc_{p}(\rho) ).$
\end{definition}

So we have defined a countable tree $\tilde{\calT} \subseteq \baire$ and, for each node $\rho \in \tilde{\calT}$, a \emph{norm} $\mu_{\rho}$ has been defined on the subsets of $\suc({\rho}).$ We can intuitively think of the norm $\mu_{\rho}$ as a way of measuring how many immediate successors has $\rho.$ Consequently, if $\mu_{ \rho}$ is big, we have more possibilities to extend subtrees at the node $\rho$. 

From the definition, it is clear that: 

\begin{lemma}\label{v50}
    For each $\rho \in \tilde{\calT}, \, \vert \suc_{\tilde{\calT}}(\rho) \vert = M(h),$ $\mu_{\rho}(\emptyset) = 0,$ $\mu_{\rho}(\suc_{\tilde{\calT}}(\rho)) = \infty$ and, if $A \subseteq \suc_{\tilde{\calT}}(\rho),$ then $\vert A \vert = \vert \suc_{\tilde{\calT}}(\rho) \vert (1 - a(h)^{-\mu_{\rho}(A)}). $
\end{lemma}

We can now define the forcing $\tilde{\bbE}$:

\begin{definition}\label{v53}
    Assume that $\tilde{\calT}, \, \varrho, \, \pi, \, a, \, M$ and $\mu_{\rho}$ are as in \autoref{v48}. The forcing $\tilde{\bbE}$ is the set $$ \tilde{\bbE} \coloneqq \left \{ p \subseteq \tilde{\calT} \colon p \text{ is a tree and } \forall \rho \in p \left( \rho \geq \trunk(p) \Rightarrow \mu_{\rho}(p) \geq 1 + \frac{1}{\lg(\trunk(p))} \right) \right \},$$ endowed with $\subseteq$. \index{$\tilde{\bbE}$}
\end{definition}


One of the important components of $\tilde{\bbE}$, apart from the $\trunk$ function, is the \emph{$\loss$ function}:

\begin{definition}\label{v57}
    Let $p \subseteq \tilde{\calT}$ be a tree. If there is a $2 < m < \omega$ such that:

    \begin{enumerate}[label=\rm{(\arabic*)}]
    
        \item $\lg(\trunk(p)) > 3m,$

        \item for any $\rho \in p, \, \lg(\rho) \geq \lg(\trunk(p)) $ entails $ \mu_{\rho}(p) \geq 1 + \frac{1}{m},$
    \end{enumerate}

    we can define $\loss(p) \coloneqq \frac{1}{m},$ where $m$ is the maximal of such $m.$ \index{$\tilde{\bbE}$!$\loss$}
\end{definition}

Therefore, $\loss$ is a function from a subset of $\tilde{\bbE}$ into $[0, 1]_{\QQQ},$ moreover (see \cite[Lem.~1.19]{KST}): 

\begin{theorem}\label{v59}\ 

    \begin{enumerate}[label=\rm{(\arabic*)}]
        \item\label{v59.1}  $\dom(\loss)$ is a dense subset of $\tilde{\bbE}.$ 
    
        \item\label{v59.2} For any $p \in \dom(\loss),$ $\frac{\Leb([p])}{\Leb([\trunk(p)])} \geq 1 - \frac{\loss(p)}{2},$ where $\Leb$ is the canonical measure on the Borel subsets of $[\tilde{\calT}].$ 
    
        \item\label{v59.3} There exists some Boolean subalgebra $\cB$ of $\calB \left(\prod_{h < \omega} M(h) \right) / \calN \left(\prod_{h < \omega} M(h) \right),$\footnote{Recall that, if $(\calX, \calA, \bfm)$ is a measure space, where $\calX$ is a Polish space and $\bfm$ is a complete measure, the collection of \emph{$\bfm$-null sets of $\calX$} is $\calN(\calX) \coloneqq \{ A \in \calA \colon \bfm(A) = 0 \}.$ } such that  $\tilde{\bbE}$ is forcing equivalent to $\cB$. More specifically, $\iota \colon \tilde{\bbE} \to \cB,$ defined by $\iota(p) \coloneqq [[p]]_{\calN},$ is a dense embedding. 
    \end{enumerate}
\end{theorem}


\section{The intersection number for forcing notions}\label{3}

In this section, based on \cite[Sec.~4.1]{uribethesis}, we introduce the notion of \emph{intersection number} for forcing notions, and we study some of its properties.  

\begin{definition}[{\cite[Def.~4.1.1]{uribethesis}}]\label{s2}
    Let $\bbP$ be a forcing notion and $Q \subseteq \bbP.$ 

    \begin{enumerate}[label=\rm{(\arabic*)}]
        \item\label{s2.1} For a finite sequence $\bar{q} = \langle q_{i} \colon i < n \rangle \in {}^{n} \bbP,$ we define $$ i_{\ast}^{\bbP}(\bar{q}) \coloneqq \max \{ \vert F \vert \colon F \subseteq n \conj  \{ q_{i} \colon i \in F\} \ \text{has a lower bound in} \ \bbP  \}. $$ 

        \item\label{s2.2} $\rm{int}^{\bbP}(Q),$ the \emph{intersection number of $Q$ in $\bbP$}, is defined by $$\Int^{\bbP}(Q) \coloneqq \inf \left \{ \frac{i_{\ast}^{\bbP}(\bar{q})}{n} \colon \bar{q} \in {}^{n} Q  \conj   n \in \omega \setminus \{ 0 \}  \right\}.$$
    \end{enumerate}

    Naturally, when the context is clear, we omit the superscript ``$\bbP$'' in $i_{\ast}^{\bbP}(\bar{q})$ and $\Int^{\bbP}(Q).$  
\end{definition}

The following result motivates the name \emph{intersection number} and it establishes a relation with the concept for Boolean algebras. 

\begin{theorem}\label{s3}
    Let $\bbP$ be a forcing notion, $\cB$ a Boolean algebra, and $\iota \colon \bbP \to \cB^{+}$ a complete embedding. If $n \in \omega \setminus \{ 0\}$ and $\bar{q} =  \langle q_{i} \colon i < n \rangle \in {}^{n} \bbP,$ then $$i_{\ast}(\bar{q}) = \max \left \{ \vert I \vert \colon I \subseteq n \conj  \bigwedge_{i \in I} \iota(q_{i}) \neq 0_{\cB} \right\}.$$
\end{theorem}

\begin{PROOF}[\textbf{Proof}]{\autoref{s3}}
    Assume that $\cB$ is a Boolean algebra, $\iota \colon \bbP \to \cB^{+}$ is a complete embedding, $n \in \omega \setminus \{ 0 \}$ and define $$m \coloneqq \max \left \{ \vert I \vert \colon I \subseteq n \conj  \bigwedge_{i \in I} \iota(q_{i}) \neq 0_{\cB} \right\}.$$ 
    
    On the  one hand, let $I \subseteq n$ be such that $b^{\ast} \coloneqq \bigwedge_{i \in I} \iota(q_{i}) \neq 0_{\cB}$ and $\vert I \vert = m.$ Since $b^{\ast} \neq 0_{\cB}$ and $\iota$ is a complete embedding, there exists $r \in \bbP$ such that it is a reduction of $b^{\ast},$ hence it is clear that $r \leq^{\bullet} q_{i}$ for all $i \in I.$ By \autoref{n5}, we can find a lower bound $q \in \bbP$ of $\{ q_{i} \colon i \in I \}$ and, therefore, $\vert I \vert \leq i^{\bbP}_{\ast}(\bar{q}).$ Thus $m \leq i_{\ast}^{\bbP}(\bar{q}).$ On the other hand, consider $F \subseteq n$ such that $\{ q_{i} \colon i \in F \}$ has a lower bound in $\bbP$, and $\vert F \vert = i_{\ast}^{\bbP}(\bar{q}).$ It is clear that $\bigwedge_{i \in I} \iota(q_{i}) \neq 0_{\cB}$ because $\{ q_{i} \colon i \in F \}$ has a lower bound and $\iota$ is a complete embedding. Thus $i_{\ast}^{\bbP}(\bar{q}) = \vert F \vert \leq m.$   
\end{PROOF}

Since  every Boolean algebra is isomorphic to a Boolean sub-algebra of $\calP(X)$ for some set $X,$ \autoref{s3} implies that  $i_{\ast}^{\cB}(\langle b_{0}, \dots, b_{n-1} \rangle)$ is the maximum number of members $b_{0}, \dots, b_{m}$  with non-empty intersection, which as mentioned above motivates the name of \emph{intersection number}. In fact, this was John Kelley's original definition in \cite{Kelley59}. 

Below we present some basic properties of the intersection number for forcing notions.

\begin{lemma}\label{s8}
    Let $\bbP$ be a forcing notion, $Q \subseteq \bbP$ such that $Q \neq \emptyset$ and $n \in \omega \setminus \{ 0 \}.$ Then,

    \begin{enumerate}[label=\rm(\arabic*)]
        \item\label{s8.1} If $Q$ is $m$-linked for some $m < \omega$  then, for all $\bar{q} \in {}^{n} Q,$ $i_{\ast}(\bar{q}) \geq \min \{ m, n \}.$ 

        \item\label{s8.2} If $\bar{q} \in {}^{n} Q,$ then $1 \leq i_{\ast}(\bar{q}) \leq n.$  As a consequence, $\Int(Q)$ is a real number and it belongs to $[0, 1].$  

        \item\label{s8.3} If  $Q$ is centered, then  for all $\bar{q} \in {}^{n} Q, \ i_{\ast}(\bar{q}) = n.$ As a consequence $\Int(Q) = 1$ if, and only if, $Q$ is centered.  

        \item\label{s8.4} For any $p \in \bbP,$ $\Int(\{ p \}) = 1.$ 
        
        \item\label{s8.5} If $Q$ is finite, then $\Int(Q) \geq \frac{1}{\vert Q \vert} > 0.$ 

        \item\label{s8.6} Let $Q$ be an anti-chain in $\bbP,$ then: 
        
        \begin{enumerate}
            \item If $Q$ is finite then $\Int(Q) = \frac{1}{\vert Q \vert}.$ 
            
            \item If $Q$ is infinite then $\Int(Q) = 0.$
        \end{enumerate}

        \item\label{s8.7} If $m \in (1, \omega)$ and $\Int(Q) \geq 1 - \frac{1}{m + 1},$ then $Q$ is $m$-linked. 

        \item\label{s8.8} If $Q \subseteq P \subseteq \bbP,$ then $\Int(P) \leq \Int(Q).$ 
    \end{enumerate}
\end{lemma}

\begin{PROOF}[\textbf{Proof}]{\autoref{s8}}\ 
    Let $n \in \omega \setminus \{ 0 \}.$ Notice that~\ref{s8.1} is clear by \autoref{s2}; \ref{s8.2} follow from~\ref{s8.1}; and~\ref{s8.4}  is a direct consequence of~\ref{s8.3}. For the other items, we must work a little bit more: 

    \begin{enumerate}
        \item[(3)] By \ref{s8.1} and \ref{s8.2} it follows that, if $Q$ is centered, then for all $\bar{q} \in {}^{n} Q,$ $i_{\ast}(\bar{q}) = n.$ Consequently,  if $Q$ is centered clearly $\Int(Q) = 1.$ On the other hand, assume that $\Int(Q) = 1$ and let $m < \omega$. Consider $\bar{q} = \langle q_{i} \colon i < m \rangle \in {}^{m} Q.$ As a consequence, $i_{\ast}(\bar{q}) \geq m,$ and therefore  $\{ q_{i} \colon i < m  \}$ has a lower bound in $\bbP$. Thus, $Q$ is $m$-linked for any $m < \omega$ and then it is centered.  
    
        \item[(5)] Assume that $Q$ is finite. Let $\bar{q} \in {}^{n} Q$ and $\ran(\bar{q}) = \{ t_{i} \colon i < k \} \subseteq Q.$ So $ k \leq \vert Q \vert.$ For each $j < k$ define $r_{j} \coloneqq \vert \{ i < n \colon q_{i} = t_{j} \} \vert,$ that is, $r_{j}$ indicates the number of times that $t_{j}$ is repeated in $\bar{q}.$ Notice that $\sum_{j < k} r_{j} = n.$ Now consider $r \coloneqq \max \{ r_{j} \colon j < k\}.$ It is clear that $ n \leq rk, \ i_{\ast}(\bar{q}) \geq r$ and, therefore, we have that $$ \frac{i_{\ast}(\bar{q})}{n} \geq \frac{r}{n} \geq \frac{r}{rk} = \frac{1}{k} \geq \frac{1}{\vert Q \vert}. $$ Thus, $\Int(Q) \geq \frac{1}{\vert Q \vert}.$ 

        \item[(6)] Let $Q$ an anti-chain in $\bbP$. 
        
        \begin{enumerate}
            \item Assume that $Q$ is finite. So $0 < \vert Q \vert < \omega.$ By the previous result it is enough to prove that $\Int(Q) \leq \frac{1}{\vert Q \vert}.$ For this, let $\bar{q} \in {}^{\vert Q \vert} Q$ be such that it lists all members of $Q.$ Since $Q$ is an anti-chain is clear that $i_{\ast}(\bar{q}) = 1.$ Therefore, $\Int(Q) \leq \frac{i_{\ast}(\bar{q})}{\vert Q \vert} = \frac{1}{\vert Q \vert}.$ Thus, $\Int(Q) \leq \frac{1}{\vert Q \vert}.$ 

            \item Assume that $Q$ is infinite. It is clear that for $m \in \omega \setminus \{0 \}$ and $\bar{q} \in {}^{m} Q$ without repetitions, $$ \Int(Q) \leq \frac{i_{\ast}(\bar{q})}{m} \leq \frac{1}{m},$$ hence $\Int(Q) \leq \frac{1}{m}$ for all $0 < m < \omega.$ Thus, by~\ref{s8.2} it follows that $\Int(Q) = 0.$
        \end{enumerate}

        \item[(7)] Suppose that $\Int(Q) \geq 1 - \frac{1}{m +1} $ for some $m \in (1, \omega).$ Towards contradiction assume that $Q$ is not $m$-linked, hence there exists a set $A = \{a_{i} \colon i < k \} \subseteq Q$ which has no lower bound and $1 < k \leq m.$ Define $\bar{q} \coloneqq \langle a_{i} \colon i < k \rangle \in {}^{k} Q.$ It is clear that $i_{\ast}(\bar{q}) < k.$ Therefore, we have that: $$1 - \frac{1}{m +1} \leq \Int(Q) \leq \frac{i_{\ast}(\bar{q})}{k} \leq \frac{k-1}{k} = 1 - \frac{1}{k}.$$
        
        This implies that $m +1 \leq k,$  which is a contradiction because $k \leq m$. Thus $Q$ is $m$-linked. 

        \item[(8)] It is clear that $Q \subseteq P$ implies $$\left \{ \frac{i_{\ast}^{\bbP}(\bar{q})}{n} \colon \bar{q} \in {}^{n} Q \conj  n \in \omega \setminus \{ 0 \}  \right\} \subseteq \left \{ \frac{i_{\ast}^{\bbP}(\bar{p})}{n} \colon \bar{p} \in P^{n} \conj  n \in \omega \setminus \{ 0 \}  \right\}.$$ Thus, by basic properties of $\inf,$ we get $\Int(P) \leq \Int(Q).$  \qedhere
    \end{enumerate}
\end{PROOF}

The characterization of $i_{\ast}^{\bbP}$ in \autoref{s3} also establishes a relation between $\Int^{\cB}(\iota(Q))$ and $\Int^{\bbP}(Q)$ for $Q \subseteq \bbP$:

\begin{corollary}\label{s4}
     Let $\bbP$ be a forcing notion, $\cB$ a Boolean algebra, $\iota \colon \bbP \to \cB^{+}$ a complete embedding and $Q \subseteq \bbP.$ If $\bar{q} = \langle q_{i} \colon i < n \rangle \in {}^{n} Q$ for $n \in \omega \setminus \{ 0 \}$ and $\bar{b} \coloneqq \langle \iota(q_{i}) \colon i < n \rangle,$   then $i_{\ast}^{\bbP}(\bar{q}) = i_{\ast}^{\cB}(\bar{b}).$ As a consequence, $\Int^{\bbP}(Q) = \Int^{\cB}(\iota(Q)).$ 
\end{corollary}

\begin{PROOF}{\autoref{s4}}
    On the one hand, let $F \subseteq n$ be a witness of $i_{\ast}^{\bbP}(\bar{q}),$ hence by \autoref{s3}, we have that $\vert F \vert = \max \left \{ \vert G \vert \colon G \subseteq n \conj  \bigwedge_{i \in G} b_{i} \neq 0_{\cB} \right\}.$ Let $I \subseteq n$ be a witness of $i_{\ast}^{\cB}(\bar{b})$, hence $\bigwedge_{i \in I} b_{i} \neq 0_{\cB}$  and therefore, $\iota_{\ast}^{\cB}(\bar{b}) = \vert I \vert \leq \vert F \vert = i_{\ast}^{\bbP}(\bar{q}).$ On the other hand, let $F$ be a witness of  $i_{\ast}^{\bbP}(\bar{q}).$ Since $\iota$ is a complete embedding, we have that $\{ b_{i} \colon i \in F \}$ has a lower bound in $\cB.$ As a consequence,  $i_{\ast}^{\cB}(\bar{b}) \leq \vert F \vert = i_{\ast}^{\bbP}(\bar{q}).$ Thus, $i_{\ast}^{\bbP}(\bar{q}) = i_{\ast}^{\cB}(\bar{b}).$ 
\end{PROOF}

As the completion of a complete Boolean algebra is the identity, we can apply  \autoref{s4} to get that the intersection number is preserved by complete embeddings for Boolean algebras: 

\begin{corollary}\label{s5}
    Let $\cB, \cC$ be complete Boolean algebras and $f \colon \cB^{+} \to \cC^{+}$ a complete embedding. If $Q \subseteq \cB,$ then $\Int^{\cB}(Q) = \Int^{\cC}(f[Q]).$ 
\end{corollary}

As a consequence, we get the analogous result for forcing notions and complete embeddings: 

\begin{corollary}\label{s6}
    Let $\bbP, \bbQ$ be forcing notions, $\iota \colon \bbP \to \bbQ$ a complete embedding and $Q \subseteq \bbP.$ Then $\Int^{\bbP}(Q) = \Int^{\bbQ}(\iota[Q]).$ As a consequence, for $R \subseteq \bbQ, \, \Int^{\bbQ}(R) \leq \Int^{\bbP}(\iota^{-1}[R]).$  
\end{corollary}

\begin{PROOF}[\textbf{Proof}]{\autoref{s6}}
    Let $\bbP, \bbQ$ be forcing notions,  $(\cB_{\bbP}, \iota_{\bbP}), (\cB_{\bbQ}, \iota_{\bbQ})$ their forcing completions, respectively, and $\iota \colon \bbP \to \bbQ$ a complete embedding. Therefore, there exists a Boolean complete embedding $f \colon \cB_{\bbP} \to \cB_{\bbQ}$ such that the following diagram commutes: 
    \[\begin{tikzcd}
    \bbP \arrow{r}{\iota} \arrow[swap]{d}{\iota_{\bbP}} & \bbQ \arrow{d}{\iota_{\bbQ}} \\
     \cB_{\bbP} \arrow{r}{f} & \cB_{\bbQ}
    \end{tikzcd}
    \]
    
    Therefore, $f \circ \iota_{\bbP} = \iota_{\bbQ} \circ \iota$ and, by applying \autoref{s4} and \autoref{s5}, we get: $$ \Int^{\bbP}(Q) = \Int^{\cB_{\bbP}}(\iota_{\bbP}[Q]) = \Int^{\cB_{\bbQ}}(f[\iota_{\bbP}[Q]]) = \Int^{\cB_{\bbQ}}(\iota_{\bbQ}[\iota[Q]]) = \Int^{\bbQ}(\iota[Q]).$$  
    
    Finally, let $R \subseteq \bbQ$. Then, $ \Int^{\bbQ}(R) \leq \Int^{\bbQ}(\iota[ \iota^{-1}[R] ]) = \Int^{\bbP}(\iota^{-1}[R]). $
\end{PROOF}

Kelley proved (see \cite[Prop.~1]{Kelley59}) that finitely additive measures can be used to define subsets of Boolean algebras whose intersection number is bounded below by a given value. 

\begin{theorem}\label{s9}
    Let $\cB$ be a Boolean algebra and   $\Xi \colon \cB \to [0, 1]$ a finitely additive measure. Consider $\bbP \coloneqq \cB^{+}$ and $\delta \in [0, 1].$ If $Q \coloneqq \{ p \in \bbP \colon \Xi(p) \geq \delta \},$ then $\Int^{\bbP}(Q) \geq \delta.$
\end{theorem}

\begin{PROOF}[\textbf{Proof}]{\autoref{s9}}
    Without loss of generality, we can assume that there exists a set $X$ such that $\cB$ is a Boolena subalgebra of $\calP(X).$ Let $n \in \omega \setminus \{ 0 \}$ and $\bar{q} = \langle q_{i} \colon i < n \rangle \in {}^{n} Q.$ For each $i < n$ consider $\chi_{i} \colon X \to \{ 0, 1\}$ as the characteristic function of $q_{i}$ in $X.$ Now, let $x \in X$ and $F_{x} \coloneqq \{ i < n \colon \chi_{i}(x) = 1 \}.$ So $\vert F_{x} \vert = \sum_{i = 0}^{n -1} \chi_{i}(x),$ $x \in \bigcap_{i \in F_{x}} q_{i}$ and $\bigcap_{i \in F_{x}} q_{i}$ is a lower bound of $\{ q_{i} \colon i \in F_{x} \}$ in $\cB^{+}.$ Hence,  by \autoref{s2}~\ref{s2.1}, we have that: $$ \sum_{i = 0}^{n -1} \chi_{i}(x) = \vert F_{x} \vert \leq i_{\ast}(\bar{q}). $$ 
    
    Therefore, for all $x \in X,$ $\sum_{i = 0}^{n -1} \chi_{i}(x) \leq i_{\ast}(\bar{q}).$  Using that $\ran(\Xi) \subseteq [0, 1],$ since by \autoref{t47}~\ref{t63} each $\chi_{i}$ is $\Xi$-integrable, we can apply the basic integral properties from \autoref{t47} to get $$i_{\ast}(\bar{q}) \geq i_{\ast}(\bar{q}) \Xi(X) =  \int_{X} i_{\ast}(\bar{q}) d \Xi \geq \int_{X} \left( \sum_{i = 0}^{n -1} \chi_{i}  \right) d \Xi  = \sum_{i = 0}^{n -1} \left( \int_{X} \chi_{i} d \Xi \right) = \sum_{i = 0}^{n-1} \Xi(q_i{})  \geq n \delta. $$ Hence $\frac{i_{\ast}(\bar{q})}{n} \geq  \delta.$ Thus $\Int(Q) \geq \delta.$ 
\end{PROOF}



Upwards closing a subset of a forcing notion does not affect its intersection number: 

\begin{lemma}\label{s7}
    If $\bbP$ is a forcing notion and $Q \subseteq \bbP,$ then $\Int^{\bbP}(Q) = \Int^{\bbP}(Q^{\uparrow}).$
\end{lemma}

\begin{PROOF}[\textbf{Proof}]{\autoref{s7}}
    Since $Q \subseteq Q^{\uparrow},$  $\Int^{\bbP}(Q^{\uparrow}) \leq \Int^{\bbP}(Q)$ follows by applying~\autoref{s8}~\ref{s8.8}. To prove the converse, let $n \in \omega \menos \{ 0 \}$ and $\bar{q} = \langle q_{i} \colon i < n \rangle \in {}^{n}(Q^{\uparrow})$. Hence, for any $i < n$, there exists $p_{i} \in Q$ such that $p_{i} \leq q_{i}.$ Then $i_{\ast}^{\bbP}(\bar{p}) \leq i_{\ast}^{\bbP}(\bar{q})$, where $\bar{p} \coloneqq \langle p_{i} \colon i < n \rangle \in {}^{n} Q$, which implies that $$ \Int^{\bbP}(Q) \leq \frac{i_{\ast}^{\bbP}(\bar{p})}{n} \leq \frac{i_{\ast}^{\bbP}(\bar{q})}{n}.$$ Finally, as $\bar{q}$ was arbitrary, it follows that $\Int^{\bbP}(Q) \leq \Int^{\bbP}(Q^{\uparrow}).$  
\end{PROOF}

\section{$\mu$-intersection-linkedness}\label{4}

Inspired in the notion of $\mu$-$\FAM$-linkedness (see \cite[Def.~4.2.8]{uribethesis}, we define: 

\begin{definition}\label{i0}
    Let $\bbP$ be a forcing notion and $\mu$ an infinite cardinal. We say that $\bbP$ is \emph{$\mu$-intersection-linked}, witnessed by the sequence $\langle Q_{\alpha, \varp} \colon \alpha < \mu \conj \varp \in (0, 1)_{\bbQ} \rangle$, if it satisfies the following conditions:  
    
    \begin{enumerate}[label=\rm{(\arabic*)}]
        \item\label{i0.1} for any $\alpha < \mu$ and $\varp \in (0, 1)_{\QQQ},$ $\Int^{\bbP}(Q_{\alpha, \varp}) \geq 1- \varp,$

        \item\label{i0.2} for any $\varp \in (0, 1)_{\QQQ},$ $\bigcup_{\alpha < \mu} Q_{\alpha, \varp}$ is dense in $\bbP.$
    \end{enumerate}
    
    If $\mu = \aleph_{0}$, we say ``$\sigma$-intersection-linked'' instead ``$\aleph_{0}$-intersection-linked''. 
\end{definition}

\begin{example}\

    \begin{enumerate}[label=\rm(\arabic*)]
        \item As a consequence of \autoref{s8}~\ref{s8.4}, any forcing notion $\bbP$ is $\vert \bbP \vert$-intersection-linked witnessed by the singletons. 

        \item Every $\mu$-$\FAM$-linked forcing notion is $\mu$-intersection-linked.  

        \item Any $\mu$-centered forcing notion is $\mu$-intersection-linked. Indeed, if $\bbP$ is $\mu$-centered, then by \autoref{s8}~\ref{s8.3} $\langle Q_{\alpha, \varp} \colon \alpha < \mu \conj \varp \in (0, 1)_{\QQQ} \rangle$ witnesses that $\bbP$ is $\mu$-intersection-linked, where for any $\alpha < \mu$ and $\varp \in (0, 1)_{\QQQ},$ $Q_{\alpha, \varp} \coloneqq Q_{\alpha}$ and $\langle Q_{\alpha} \colon \alpha < \mu \rangle$ witnesses that $\bbP$ is $\mu$-centered. In particular, Hechler forcing $\bbD,$ the canonical forcing notion to add dominating reals, is $\sigma$-intersection-linked. As a consequence, by \autoref{n80}, it is an example of a $\sigma$-intersection-linked forcing notion that is no $\sigma$-$\FAM$-linked because it is no $\sigma$-$\Fr$-linked (see \cite[Cor.~4.2.11]{uribethesis}). 
    \end{enumerate}
\end{example}

The following observation shows that we can change condition \ref{i0.2} in \autoref{i0}: 

\begin{observation}\label{i3}
     Notice that, by~\autoref{s7}, a forcing notion $\bbP$ is $\mu$-intersection-linked witnessed by the sequence $\langle Q_{\alpha, \varp} \colon \alpha < \mu \conj \varp \in (0, 1)_{\QQQ} \rangle$ if, and only if, it is witnessed by the sequence $\langle Q_{\alpha, \varp}^{\uparrow} \colon \alpha < \mu \conj \varp \in (0, 1)_{\QQQ} \rangle,$ since $\bigcup_{\alpha < \mu} Q_{\alpha, \varp}$ is dense in $\bbP$ if, and only if, $\bigcup_{\alpha < \mu} Q_{\alpha, \varp}^{\uparrow} = \bbP,$ for all $\varp \in (0, 1)_{\QQQ}.$ Consequently, we have an equivalent statement of \autoref{i0}: we can replace condition \autoref{i0} \ref{i0.2} by ``for every $\varp \in (0, 1)_{\QQQ},$ $\bigcup_{\alpha < \mu} Q_{\alpha, \varp} = \bbP$''. That is, if $\bbP$ is $\mu$-intersection-linked witnessed by the sequence $\langle Q_{\alpha, \varp} \colon \alpha < \mu \conj \varp \in (0, 1)_{\QQQ} \rangle,$ we can assume without loss of generality,  that each $Q_{\alpha, \varp}$ is upwards closed, hence \autoref{i0} \ref{i0.2} can be changed by ``for every $\varp \in (0, 1)_{\QQQ}, $ $\bigcup_{\alpha < \mu} Q_{\alpha, \varp} = \bbP$''.  
\end{observation}

Recall that $\varphi(x)$ is \emph{hereditary forcing property} if for any pair of forcing notions $\bbP,$ $\bbQ$ such that $\bbP \lessdot \bbQ,$ if $\varphi(\bbQ)$ holds, then $\varphi(\bbP)$ holds. 

\begin{theorem}\label{i7}
    $\mu$-intersection-linkedness is a hereditary forcing property. 
\end{theorem}

\begin{PROOF}{\autoref{i7}}
    Let $\iota \colon \bbP \to \bbQ$ be a function. On the one hand, assume that $\iota$ is a dense embedding and that $\bbP$ is $\mu$-intersection-linked witnessed by $\langle P_{\alpha, \varp} \colon \alpha < \mu \conj \varp \in (0, 1)_{\QQQ} \rangle.$ For any $\alpha < \mu$ and $\varp \in (0, 1)_{\QQQ},$ define $Q_{\alpha, \varp} \coloneqq \iota[P_{\alpha, \varp}]$ and let us show that $\langle Q_{\alpha, \varp} \colon \alpha < \mu \conj \varp \in (0, 1)_{\QQQ} \rangle$ witnesses that $\bbQ$ is $\mu$-intersection-linked. For the first condition, fix $\alpha < \mu$ and $\varp \in (0, 1)_{\QQQ}.$ By \autoref{s6}, we have that $1 - \varp \leq \Int^{\bbP}(P_{\alpha, \varp}) = \Int^{\bbQ}(\iota[P_{\alpha, \varp}]) = \Int^{\bbQ}(Q_{\alpha, \varp}).$ For the second condition, let $\varp \in (0, 1)_{\bbQ}$ and $q \in \bbQ.$ By density, we can find $p \in \bbP$ such that $\iota(p) \leq q,$ and by \autoref{i0}~\ref{i0.1}, there are $\beta < \mu$ and $p' \in P_{\beta, \varp}$ such that $p' \leq p$, therefore $\iota(p') \leq \iota(p) \leq q$ and $\iota(p') \in Q_{\beta, \varp},$ hence $\bigcup_{\alpha < \mu} Q_{\alpha, \varp}$ is dense in $\bbQ$. Thus, $\bbQ$ is $\mu$-intersection-linked. On the other hand, assume that $\iota$ is a complete embedding and that $\bbQ$ is $\mu$-intersection-linked witnessed by $\langle Q_{\alpha, \varp} \colon \alpha < \mu \conj \varp \in (0, 1)_{\QQQ} \rangle.$ For any $\alpha < \mu$ and $\varp \in (0, 1)_{\QQQ},$ define $P_{\alpha, \varp} \coloneqq \iota^{-1}[Q_{\alpha, \varp}].$ Let us show that $\langle P_{\alpha, \varp} \colon \alpha < \mu \conj \varp \in (0, 1)_{\QQQ} \rangle$ witnesses that $\bbP$ is $\mu$-intersection-linked. For the first condition, fix $\alpha < \mu$ and $\varp \in (0, 1)_{\QQQ}.$ By \autoref{s6} we have that, $1 - \varp \leq \Int^{\bbQ}(Q_{\alpha, \varp}) \leq \Int^{\bbP}(\iota^{-1}[Q_{\alpha, \varp}]) = \Int^{\bbP}(P_{\alpha, \varp}).$ For the second condition, let $p \in \bbP$ and $\varp > 0.$ By \autoref{i3}, we can find $\beta < \mu$ such that $\iota(p) \in Q_{\beta, \varp}$, hence $p \in P_{\beta, \delta}.$ Thus, $\bigcup_{\alpha < \mu} P_{\alpha, \varp} = \bbP.$ Finally, $\bbP$ is $\mu$-intersection-linked. 
\end{PROOF}

Recall that $\varphi(x)$ is a \emph{forcing property} if for any pair of forcing-equivalent forcing notions $\bbP,$ $\bbQ,$ we have that $\varphi(\bbP)$ if, and only if, $\varphi(\bbQ)$. Notice that, any hereditary property is a forcing property. Consequently, 

\begin{corollary}\label{i8}
    $\mu$-intersection-linkedness is a forcing property. 
\end{corollary}

The property of $\mu$-intersection-linkedness is stronger than $\mu$-$m$-linked for any $m < \omega$: 

\begin{theorem}\label{i15}
    Every $\mu$-intersection-linked forcing notion is $\mu$-$m$-linked for all $0 < m < \omega.$ 
\end{theorem}

\begin{PROOF}[\textbf{Proof}]{\autoref{i15}}
    Let $\bbP$ be a $\mu$-intersection-linked forcing notion witnessed by $\langle Q_{\alpha, \varp} \colon \alpha < \mu \conj \varp \in (0, 1)_{\QQQ} \rangle,$ and $m < \omega$ such that $m > 0.$ Since $ \frac{1}{m+1} \in (0, 1),$ there exists some $\varp_{m} \in (0, 1)_{\QQQ}$ such that $\varp_{m} < \frac{1}{m+1}.$ Therefore, by~\autoref{i0}~\ref{i0.1}, for all $\alpha < \mu$ we have that, $$ 1 - \frac{1}{m+1} < 1 - \varp_{m} \leq \Int(Q_{\alpha, \varp_{m}}). $$
    
    Applying  \autoref{s8}~\ref{s8.7}, we get that for all $\alpha < \mu,$ $Q_{\alpha, \varp_{m}}$ is $m$-linked. It is clear that each $Q_{\alpha, \varp_{m}}^{\uparrow}$ is $m$-linked because $Q_{\alpha, \varp_{m}}$ is. Finally, let $p \in \bbP.$ By \autoref{i0}~\ref{i0.2}, there are $\alpha_{0} < \mu$ and $q \in Q_{\alpha_{0}, \varp_{m}}$ such that $q \leq p,$ so $p \in Q_{\alpha_{0}, \varp_{m}}^{\uparrow},$ hence $\bigcup_{\alpha < \mu} Q_{\alpha, \varp_{m}}^{\uparrow} = \bbP.$ Thus, $\bbP$ is $\mu$-$m$-linked.     
\end{PROOF}

\begin{corollary}\label{i15.1}
    Every $\mu$-intersection-linked forcing notion is $\theta$-$m$-Knaster for all $0 < m < \omega$ and any cardinal $\theta$ with $\cf(\theta)>\mu$. 
\end{corollary}

Recall that, for $H \subseteq \baire,$ we can define the \emph{localization relational system} $\Lc(\omega,H) \coloneqq \la\omega^\omega, \Swf(\omega,H), \in^*\ra,$ where for $h \in \baire$ the sets of \emph{slaloms} $\calS(\omega, h)$, $\calS(\omega, H)$ are defined by: $$ \Swf(\omega,h)  \coloneqq \prod_{i<\omega}[\omega]^{\leq h(i)} \text{ and } \Swf(\omega,H) \coloneqq \bigcup_{h\in H}\Swf(\omega,h),$$ and  for functions $x, y \in \baire$, the relation ``\emph{$y$ localizes $x$}'' is giving by: $x\in^* y$ iff there exists $ m<\omega $ such that for any $i\geq m $, we have $x(i)\in y(i).$

Inspired by a result of Kamburelis~\cite{Ka}, in an unpublished result, Diego A. Mej\'ia showed that, if a forcing notion can be covered with $\theta$ subsets with positive intersection number, then it is $\mu^+$-$\Lc(\omega, H_*)$-good, where $H_* \coloneqq \set{\id^{k+1}}{k<\omega}$.\footnote{Powers of the identity function on $\omega$.} As a consequence,  

\begin{theorem}
    Any $\mu$-intersection-linked forcing notion is $\mu^+$-$\Lc(\omega,H_*)$-good. 
\end{theorem}

\begin{corollary}
    Any $\mu$-centered forcing notion is $\mu^+$-$\Lc(\omega,H_*)$-good.
\end{corollary}

We close this paper by showing some non-trivial examples of $\mu$-intersection-linked forcing notions. First, inspired in \cite{MU23}, we have: 

\begin{theorem}\label{i70}
    For any infinite cardinal $\theta,$ any Boolean algebra with a strictly positive probability measure satisfying the $\theta$-density property is $\theta$-intersection-linked.
\end{theorem}

\begin{PROOF}{\autoref{i70}}
     Let $\cB$ a Boolean algebra, $\mu$ a strictly positive probability measure on it, and fix $S \subseteq \cB^{+}$ witnessing that $\mu$ satisfies the $\theta$-density property. For each  $s \in S$ and $\varp \in (0, 1)_{\bbQ},$ consider the sets $Q_{s, \varp} \coloneqq \{ b \in \cB \colon \mu_{s}(b) \geq 1 - \varp \}.$ Let us prove that the sequence $\langle Q_{s, \varp} \colon s \in S \conj  \varp \in (0, 1)_{\bbQ} \rangle$ witnesses that $\cB$ is $\theta$-intersection-linked by proving the conditions in \autoref{i0}. For the first one, notice that by \autoref{s9}, for any $\varp \in (0, 1)_{\QQQ}$ and $s \in S,$ we have that $\Int(Q_{s, \varp}) \geq 1 - \varp.$ For the second one, let $\varp \in (0, 1)_{\bbQ}$ and $b \in \cB^{+}.$ By the $\theta$-density property, there exists some $s \in S$ such that, $\mu_{s}(b) \geq 1 - \varp,$ that is, $b \in Q_{s, \varp}.$ Thus, not only $\bigcup_{s \in S} Q_{s, \varp_{0}}$ is dense in $\cB$ but it is the whole $\cB$. Thus, $\cB$ is $\mu$-intersection-linked. 
\end{PROOF}  

Consequently, by \autoref{i3.0}, we get: 

\begin{corollary}\label{i99}
    For any infinite cardinal $\theta,$ the random algebra adding $\theta$-many random reals is $\theta$-intersection-linked. In particular, random forcing is $\sigma$-intersection-linked. 
\end{corollary}

\begin{theorem}\label{i200}
    $\tilde{\bbE}$ is $\sigma$-intersection-linked.
\end{theorem}

\begin{PROOF}{\autoref{i200}}
    For any $t \in \tilde{\calT}$ and $\varp \in (0, 1)_{\QQQ}$, consider the sets $$Q_{t, \varp} \coloneqq \{ p \in \dom(\loss) \colon \trunk(p) = t \conj \loss(p) \leq \varp \}.$$ Since $\tilde{\calT}$ is countable, it is enough to prove that $\langle Q_{t, \varp} \colon t \in \tilde{\calT} \conj \varp \in (0, 1)_{\QQQ} \rangle$ satisfies conditions \ref{i0.1} and~\ref{i0.2} in~\autoref{i0}. 

    On the one hand, for condition \ref{i0.1}, fix $t \in \tilde{\calT}$ and $\varp \in (0, 1)_{\bbQ}.$ By \autoref{v59}~\ref{v59.3}, $\iota \colon \tilde{\bbE} \to \cB$  such that for all $p \in \tilde{\bbE},$ $\iota(p) \coloneqq [[p]]_{\calN}$ is a dense embedding. Also, by \autoref{v59}~\ref{v59.2}, for all $p \in \dom(\loss),$ we have that: $$ \frac{\Leb([p])}{\Leb([\trunk(p)])} \geq 1 - \frac{\loss(p)}{2}.$$

    In particular, if $p \in Q_{t, \varp}$ then $\trunk(p) = t, \loss(p)  \leq \varp$ and therefore, for any $p \in Q_{t, \varp},$ $  \Leb_{[t]}([p]) \geq 1 - \varp$. Define $\mu \colon \cB \to [0, 1]$ such that for any $b \in \cB,$ $\mu([b]_{\calN}) \coloneqq \Leb_{[t]}(b).$ Then, $\mu([\emptyset]_{\calN}) = \Leb_{[t]}(\emptyset) = 0.$ Also, if $[a]_{\calN}, [b]_{\calN} \in \cB$  are such that $[a]_{\calN} \conj  [b]_{\calN} = [\emptyset]_{\calN},$ then $[a \cap b]_{\calN} = [\emptyset]_{\calN},$ hence $ a \cap b = (a \cap b) \triangle \emptyset \in \calN.$ Therefore, 
        \begin{equation*}
            \begin{split}
                \mu([a]_{\calN} \vee [b]_{\calN}) &= \mu([a \cup b]_{\calN})  = \Leb_{[t]}(a) + \Leb_{[t]}(b) - \Leb_{[t]}(a \cap b)\\
                & = \Leb_{[t]}(a) + \Leb_{[t]}(b) - \Leb_{[t]}(\emptyset) =  \Leb_{[t]}(a) + \Leb_{[t]}(b)\\
                & = \mu_{[t]}([a]_{\calN}) + \mu_{[t]}([b]_{\calN}).
            \end{split}
        \end{equation*}

    Thus, $\mu$ is a finitely additive measure. 
        
    Now, define also $Q \coloneqq \{ b \in \cB \colon \mu(b) \geq 1 - \varp \}.$ Notice that $\iota[Q_{t, \varp}] \subseteq Q.$ Indeed, let $b \in \iota[Q_{t, \varp}]$ and $p \in Q_{t, \varp}$ such that $b = [p]_{\calN}.$ Since $t = \trunk(p),$ $[p] \subseteq [t]$ hence, $[p] \cap [t] = [p].$ Therefore, we have that $ \mu(b) = \mu([[p]]_{\calN}) = \Leb_{[t]}([p]) \geq 1 - \varp,$ that is, $b \in Q,$ hence $\iota[Q_{t, \varp}] \subseteq Q.$ Finally, by \autoref{s8}~\ref{s8.8}, \autoref{s9} and~\autoref{s6} we have that: $$ \Int^{\tilde{\bbE}}(Q_{t, \varp}) = \Int^{\cB}(\iota[Q_{t, \varp}]) \geq \Int^{\cB}(Q) \geq 1 - \varp.$$  Thus, $\Int^{\tilde{\bbE}}(Q_{t, \varp}) \geq 1 - \varp.$

    On the other hand, for condition \ref{i0.2}, let $\varp \in (0, 1)_{\QQQ}$ and $p \in \tilde{\bbE}.$ By \autoref{v59}~\ref{v59.1}, there is some $q \in \dom(\loss)$ such that $q \leq p$ and, by the definition of $\tilde{\bbE},$ we can extend $\trunk(q)$ long enough to find a condition $r \in \tilde{\bbE}$ such that $r \leq q$ and $\loss(r) \leq \varp.$ So if we define $t \coloneqq \trunk(r),$ we have that $r \in Q_{t, \varp}$ and $r \leq p.$
\end{PROOF}

\begin{remark}\label{i201}
    An alternative and easier proof of \autoref{i200} is to apply~\autoref{i70} and ~\autoref{m15.1}. 
\end{remark}

\section*{Acknowledgements}

This paper was developed for the conference proceedings corresponding to the \emph{RIMS Set Theory Workshop 2023: Large Cardinals and the Continuum}, organized by Professor Hiroshi Fujita, where the author participated in giving the talk titled ``\emph{A general theory of iterated forcing using finitely additive measures}''. 

Much of the content of this paper is part of the author's master's thesis. So the author is deeply grateful to his main advisor, Professor Diego A. Mejía, from Shizuoka University, for his wonderful and careful advice in the conception, development, and writing of the thesis. Also, for reviewing this document, helping to improve its presentation, and for pointing out \autoref{i201}. 

The author is especially thankful to Miguel A. Cardona and Diego A. Mej\'ia, for having pointed out the main idea in the proof of \autoref{i200} and for valuable discussions.

This work was partially supported by the Austrian Science Fund (FWF) through the project P33895.

{\small
\bibliography{appl}
\bibliographystyle{alpha}
}

\end{document}